\documentclass[11pt]{article}

\usepackage[T1]{fontenc}
\usepackage[utf8]{inputenc}

\usepackage{fourier} %pour le symbole \danger
\usepackage{soul} %pour souligner \ul avec respect des retours à la ligne
\usepackage{enumerate} %pour les \begin{enumerate}[(i)] par exemple

\usepackage{amsfonts,amssymb,amsmath,amsthm}
\usepackage{bbm}
\usepackage{pgf}

\numberwithin{equation}{section}

\def\11{\mathbbm{1}}

\def\E{\mathbb{E}}
\def\P{\mathbb{P}}
\def\R{\mathbb{R}}

\def\N{\mathbb{N}}

\def\d{\partial}
\def\Z{\mathbb{Z}}

\def\cE{{\cal E}}

\newtheorem{thm}{Theorem}[section]

\newtheorem{prop}[thm]{Proposition}

\theoremstyle{remark}
\newtheorem{rem}{Remark}

\newcommand{\vertiii}[1]{{\left\vert\kern-0.25ex\left\vert\kern-0.25ex\left\vert #1 
    \right\vert\kern-0.25ex\right\vert\kern-0.25ex\right\vert}}
\begin{document}

\title{Convergence of the Fleming-Viot process toward the minimal
  quasi-stationary distribution}

\author{Nicolas Champagnat$^{1,2,3}$, Denis Villemonais$^{1,2,3}$}

\footnotetext[1]{IECL, Universit\'e de Lorraine, Site de Nancy, B.P. 70239, F-54506 Vandœuvre-lès-Nancy Cedex, France}
\footnotetext[2]{CNRS, IECL, UMR 7502, Vand{\oe}uvre-l\`es-Nancy, F-54506, France}  
\footnotetext[2]{Inria, TOSCA team, Villers-l\`es-Nancy, F-54600, France.\\
  E-mail: Nicolas.Champagnat@inria.fr, Denis.Villemonais@univ-lorraine.fr}

\maketitle

\begin{abstract}
  We prove under mild conditions that the Fleming-Viot process selects
  the minimal quasi-stationary distribution for Markov
  processes with soft killing on non-compact state spaces. Our results
  are applied to multi-dimensional birth and death processes, 
  continuous time Galton-Watson processes and diffusion processes
  with soft killing.
\end{abstract}



\section{Introduction}
\label{sec:intro}

Let
$\left(\Omega,(\mathcal{F}_{t})_{t\in[0,+\infty)},\P,(X_t)_{t\in
    [0,+\infty]}\right)$ be a continuous time Markov process evolving
in a Polish state space $(E,\cE)$. Let also $\kappa:E\rightarrow\R_+$ be a bounded measurable
function and consider the killed Markov process
$(Y_t)_{t\in [0,+\infty)}$ evolving in $E\cup\{\d\}$, where
$\d\notin E$ is a cemetery point, defined by
\begin{align*}
Y_t=\begin{cases}
X_t&\text{ if }\int_0^t \kappa(X_s)\,ds<\zeta,\\
\d&\text{ otherwise,}
\end{cases}
\end{align*}
where $\zeta$ is an exponential random variable with parameter $1$,
independent from $X$. A quasi-stationary distribution for $Y$ is a
probability measure $\nu_{QSD}\ $ on $E$ such that
\begin{align*}
  \lim_{t\rightarrow+\infty} \P_{\mu}(Y_t\in\cdot\,\mid\,t<\tau_\d)=\nu_{QSD},
\end{align*}
for some probability measure $\mu$ on $E$, and where the convergence
holds \textcolor{black}{for the total variation distance (for example)}. The distribution $\nu_{QSD}$ is called the
\textit{minimal quasi-stationary distribution} of $Y$ if the above
convergence holds true for all $\mu=\delta_x$, where $x$ runs over
$E$~\textcolor{black}{\cite{AsselahThai2012,AsselahFerrariEtAl2016}}.

This paper is devoted to the study, in a general setting, of the
ergodicity and the convergence of a Fleming-Viot type particle system
in the long time/large particle limit toward the minimal
quasi-stationary distribution of $Y$. In particular, we recover and
sensibly improve the results of~\cite{AsselahFerrariEtAl2016}
and~\cite{Villemonais2015}, where adhoc probabilistic and spectral
considerations were used to study this problem for continuous time
Galton-Watson processes and birth and death processes
respectively. The main difficulty in these examples is that the state
space is neither compact nor the process comes back from infinity (see
also~\cite{AsselahThai2012} for a fine control of the right most
particle position in a Fleming-Viot type model). In order to control
the behaviour of the particle system despite the lack of compactness,
we make use of the general results
of~\cite{ChampagnatVillemonais2017b} for the study of quasi-stationary
distributions. Our main contributions \textcolor{black}{compared to the previously cited references} are the generality of our
approach (our framework includes, for instance, the case of
multi-dimensional birth and death processes, of multi-type Galton
Watson processes and of diffusion processes on $\R^d$),  a new 
speed of convergence result and a stronger convergence.

% satisfying an integrability condition and where the convergence holds with respect to the total variation norm.

% In this paper, we first state a sufficient criterion ensuring that the process $Y$ admits a quasi-stationary distribution, directly derived from those obtained in~\cite{ChampagnatVillemonais2017b} (in the discrete state space case, one can also use the $R$-positive matrix theory approach developed in ....). This criterion ensures that there exists a probability measure $\nu_{QSD}$ on $E$ such that
% \begin{align*}
% \lim_{t\rightarrow+\infty} \P_{\mu}(Y_t\in\cdot\,\mid\,t<\tau_\d)=\nu_{QSD},
% \end{align*}
% for all probability measure $\mu$ on $E$ satisfying an integrability condition and where the convergence holds with respect to the total variation norm.

 Let us describe informally the dynamics of the Fleming-Viot particle system with $N\geq 2$ particles, which we denote by
$(X^1,X^2,\ldots,X^N)$. The process starts at a position
$(X^1_0,X^2_0,\ldots,X^N_0)\in E^{N}$ and evolves as follows:
\begin{itemize}
\item[-] the particles $X^i$, $i=1,\ldots,N$, evolve as independent copies of the process $Y$ until one of them hits $\d$; this hitting time is denoted by $\tau_1$;
\item[-] then the (unique) particle hitting $\d$ at time $\tau_1$ jumps instantaneously on the position of a particle chosen uniformly among the $N-1$ remaining ones; this operation is called a \textit{rebirth};
\item[-]  because of this rebirth, the $N$ particles lie in $E$ at time $\tau_1$; then the $N$ particles evolve as independent copies of $Y$ and so on.
\end{itemize}

% \medskip \denis{J'ai écrit la partie qui suit pour se placer dans la
%   plus grande g\'en\'eralit\'e possible, cependant ce n'est pas très
%   satisfaisant \`a mon go\^ut... En effet, je ne vois pas comment
%   garantir que le processus ci-dessous, en supposant qu'il existe et
%   qu'il est cadlag, correspond au processus d\'ecrit informellement
%   ci-dessus... Peut-\^etre que la description du semi-groupe peut
%   aider, mais je n'en suis m\^eme pas s\^ur...}

\medskip More formally, for a fixed number of particles $N\geq 2$, we consider the Markov process in $E^N$ whose weak infinitesimal
generator \textcolor{black}{(in the sense of~\cite{MeynTweedie1993}, see also~\cite{ChampagnatVillemonais2017} for an even weaker
  version, which could also be used in this context)} can be expressed as
\begin{align*}
  L^N f(x_1,\ldots,x_N)=\sum_{i=1}^N L\varphi(x_i) + \frac{\kappa(x_i)}{N-1}\sum_{j=1,\,j\neq i}^N (\varphi(x_j)-\varphi(x_i)),
\end{align*}
where $L$ is the weak infinitesimal generator of $X$ and for all $f$
of the form $f(x_1,\ldots,x_N)=\sum_{i=1}^N \varphi(x_i)$, with
$\varphi$ belonging to the \textcolor{black}{extended} domain of $L$ \textcolor{black}{(in the sense of~\cite{MeynTweedie1993})}. We
assume that an almost 
surely c\`adl\`ag strong Markov process with such an infinitesimal
generator exists. In most cases, it is easy to build such a Markov
process, but imposing general conditions on $X$ and $\kappa$ ensuring
this is challenging. Let us simply mention that, if $X$ is a regular
non-explosive pure jump Markov process in $E$, then $L^N$ corresponds
to the infinitesimal generator of a non-explosive pure jump Markov
process in $E^N$ and the Fleming-Viot process is thus
well-defined. Also, if $X$ is the solution to a stochastic differential equation with smooth
coefficients in $E=\R^d$, then $L^N$ corresponds to the infinitesimal
generator of the solution on $\R^{d\times N}$ to a stochastic
differential equation with jumps occuring at bounded jump rate, which also admits an
almost sure c\`adl\`ag representation. % Finally, if $X$ is a Feller
% process and $\kappa$ is continuous, then $L^N$ corresponds to the
% infinitesimal generator of a Feller Markov process on $E^N$.
% (\denis{\`a v\'erifier, mais je crois que c'est faisable en
%   d\'ecrivant le semi-groupe explicitement, ou en utilisant les
%   th\'eor\`emes g\'en\'eraux \`a la Ethier Kurtz p.165 Theorem 2.2 par
%   exemple}).

\medskip This Fleming-Viot type system has been introduced by Burdzy,
Holyst, Ingermann and March in \cite{BurdzyHolystEtAl1996} and studied
in \cite{BurdzyHolystEtAl2000}, \cite{GrigorescuKang2004},
\cite{Villemonais2011}, \cite{GrigorescuKang2012} for
multi-dimensional diffusion processes. The study of this system when
the underlying Markov process $X$ is a continuous time Markov chain in
a countable state space has been initiated in \cite{FerrariMaric2007}
and followed by \cite{AsselahFerrariEtAl2011},
\cite{AsselahFerrariEtAl2016}, \cite{GroismanJonckheere2013},
\cite{AsselahThai2012} and \cite{CloezThai2016}. We also refer the
reader to \cite{GroismanJonckheere2013a}, where general considerations
on the link between the study of such systems and front propagation
problems are considered and
to~\cite{CerouDelyonEtAl2016,DelyonCerouEtAl2017} where CLTs for this
Fleming-Viot type process have been proved.

Note that our settings correspond to the \textit{soft} killing case,
so that, denoting by $\tau_1<\tau_2<\cdots<\tau_n<\cdots$ the sequence
of rebirths times and since the rate at which rebirths occur is
uniformly bounded above by $N\,\|\kappa\|_\infty$,
\begin{align*}
\lim_{n\rightarrow\infty} \tau_n=+\infty\text{ almost surely.}
\end{align*}
Hence, contrarily to the \textit{hard} killing
case~\cite{BieniekBurdzyEtAl2011}, there is no risk for the particle
system $(X^1_t,X^2_t,\ldots,X^N_t)_{t\geq 0}$ to undergo an infinite
number of jump in finite time. This guarantees that it
is well defined for any time $t\geq 0$ in an incremental way (rebirth
after rebirth).

\medskip We emphasise that, because of the rebirth mechanism, the
particle system $(X^1,X^2,\ldots,X^N)$ evolves in $E^N$. For any
$t\geq 0$, we denote by $\mu^N_t$ the empirical distribution of
$(X^1,X^2,\ldots,X^N)$ at time $t$, defined by
\begin{align*}
\mu^N_t=\frac{1}{N}\sum_{i=1}^N\delta_{X^i_t}\in{\cal M}_1(E),
\end{align*}
where ${\cal M}_1(E)$ is the set of probability measures on $E$. Our
aim is to find a tractable condition ensuring: 1) that $Y$ admits a
minimal quasi-stationary distribution using the recent general results
of~\cite{ChampagnatVillemonais2017b}; 2) that the law of $\mu^N_t$
converges toward the law of a (random) distribution ${\cal X}^N$ on
$E$ using classical Foster-Lyapunov type
criteria~\cite{MeynTweedie1993}; and 3) that this sequence of random
distributions converges to the (deterministic) minimal quasi-stationary distribution of $Y$
using the general convergence results of~\cite{Villemonais2014}.

% A
% general convergence result obtained in \cite{Villemonais2014} ensures
% that, if $\mu_0^N\rightarrow\mu_0$ with respect to the weak topology,
% then
% \begin{align*}
% \mu^N_t\xrightarrow[N\rightarrow\infty]{} \P_{\mu_0}(X_t\in\cdot\mid t<T_0).
% \end{align*}
% However, the generality of this result does not extend to the long time behaviour of the particle system, which is the subject of the next result.

%  We then use these recent results to prove the existence of the long time and high number of particles limit of a Fleming-Viot type particle system. The particles of
% this system evolve as independent copies of the killed process $Y$, but they
% undergo rebirths when they hit $\d$ instead of being trapped at the cemetery point (in particular,
% the number of particles that are in $E$ remains constant as time goes on). Our main
% result is that the empirical stationary distribution of the particle system exists and converges to the minimal quasi-stationary distribution of the
% underlying process $Y$.

Our main results are stated in Section~\ref{sec:mainresults},
illustrated by several examples in Section~\ref{sec:examples} and
proved in Section~\ref{sec:proofs}.

\section{Main results}
\label{sec:mainresults}

We first present our main assumptions. Several examples satisfying this
requirement are provided in the next section. For the definition of
the \textit{extended domain} of the generator, we refer the reader to
Meyn \& Tweedie~\cite{MeynTweedie1993}. We recall that a
subset $K\in E$ is called a \textit{small set} for $X$ if there exist a positive time
$t_K>0$, a positive constant $\alpha_K>0$ and a probability measure
$\nu_K$ on $K$ such that, for all $x\in K$,
\[
  \P_x(X_{t_K}\in A)\geq \alpha_K\,\nu_K(A\cap K),\quad\forall A\subset E.
\]
We also refer the reader to~\cite{MeynTweedie1993} and references
therein for general considerations on small sets and the related
concept of~\textit{petite sets}.

\medskip

\textbf{Assumption H.}  Assume that all compact sets are small sets
for $X$. Assume also that there exists a locally bounded function
$V:E\rightarrow[1,+\infty)$ in the extended domain of the
infinitesimal generator of $X$, with relatively compact level sets and
such that
\begin{align*}
LV(x)\leq -\lambda_1 V(x)+C,\ \forall x\in E,
\end{align*}
for some constant $\lambda_1>\|\kappa\|_\infty$.  Finally, assume that
$Y$ satisfies the following property: for all compact set
$K\subset E$,
\begin{align*}
\inf_{t\geq 0}\frac{\inf_{x\in K}\P_x\left(Y_t\notin\d\right)}{\sup_{x\in K}\P_x\left(Y_t\notin\d\right)}>0.
\end{align*}

\bigskip

Of course, the two first points of Assumption~H imply that $X$
satisfies a Foster-Lyapunov type criterion and
hence that it is exponentially ergodic. However, this does not
guarantee the long-time convergence of the conditional distribution of
$Y$. For instance, if $X$ is a birth and death process evolving in
$\N=\{1,2,3,\ldots\}$ and if $\kappa(x)=\11_{x=1}$, the exponential
ergodicity of $X$ is not sufficient to deduce the long-time
convergence of the conditional distribution of $Y$ (as clearly appears
in the reference work~\cite{DoornErik1991}). The last point of
Assumption~H allows us to overcome this difficulty (using the recent
results of~\cite{ChampagnatVillemonais2017b}).

\begin{thm}
\label{thm:main-QSD}
If Assumption~H holds true, then $Y$ admits a unique quasi-stationary
distribution $\nu_{QSD}$ such that
$\nu_{QSD}(V)<+\infty$. Moreover, there exist two positive
constants $c,\gamma>0$ such that
\begin{align*}
\left\|\P_{\mu}(Y_t\in\cdot\mid t<\tau_\d)-\nu_{QSD}\right\|_{TV}\leq c\,\mu(V)\,e^{-\gamma t},\ \forall t>0,
\end{align*}
for any probability distribution $\mu$ on $E$ such that $\mu(V)<\infty$.
\end{thm}

Since $V(x)<\infty$ for all $x\in E$, $\nu_{QSD}$ is the minimal quasi-stationary distribution of $Y$.
Several other properties can be deduced from Assumption H and
Theorem~\ref{thm:main-QSD}, as detailed
in~\cite[Section 2]{ChampagnatVillemonais2017b}.

% \begin{rem}
%   Since any measure $\mu$ with compact support satisfies
%   $\mu(V)<\infty$, we deduce that the QSD obtained in
%   Theorem~\ref{thm:main-QSD} is the so-called \textit{Yaglom limit}
%   (see~\cite{MeleardVillemonais2012}), also referred to as the
%   \textit{minimal QSD} (see~\cite{AsselahFerrariEtAl2016}).
% \end{rem}

\begin{rem}
  One can weaken Assumption~H replacing $\lambda_1>\|\kappa\|_\infty$
  by $\lambda_1>\text{osc }\kappa$, where $\text{osc }\kappa=\sup_{x\in E}\kappa(x)-\inf_{x\in E}\kappa(x)$ is the oscillation of $\kappa$. However, this generalisation does
  not transfer to the next results.
\end{rem}

The next theorem states that the Fleming-Viot process defined in the
introduction is exponentially ergodic for $N$ large enough.

\begin{thm}
\label{thm:main-ergodicity}
Assume that Assumption~H is satisfied. Then, for any
$N> \frac{\lambda_1}{\lambda_1-\|\kappa\|_\infty}$, the process
$(X^1,\ldots,X^N)_{t\geq 0}$ is exponentially ergodic, which means that there
exists a probability measure $M^N$ on $E^N$ such that
\begin{align*}
  \left\|\text{Law}(X^1_t,\ldots,X^N_t\,\mid\, (X^1_0,\ldots,X^N_0)=x\in E^N)-M^N\right\|_{TV} \leq C_N(x) e^{-\gamma_N t},
\end{align*}
where $C_N:E^N\rightarrow [0,+\infty)$ is a non-negative function and $\gamma_N>0$ is a constant.
Moreover
\[
  \int_{E^N}\sum_{i=1}^N V(x_i)\,dM^N(x_1,\ldots,x_N)\leq C/(\lambda_1-\|\kappa\|_\infty N/(N-1)).
\]
% and
% \begin{align*}
% \cX^N\xrightarrow[N\rightarrow\infty]{Law}\nu_{QSD},
% \end{align*}
% where $\nu_{QSD}$ is the minimal QSD of $Y$.
\end{thm}

Note that, using the above theorem, one deduces that $\mu^N_t$ also
converges in law, when $t\rightarrow+\infty$, to a (random)
probability measure ${\cal X}^N$ on $E$, defined as
\begin{align*}
  {\cal X}^N=\frac{1}{N}\sum_{i=1}^N \delta_{x_i},\ \text{where }\text{Law}(x_0,\ldots,x_N)=M^N.
\end{align*}

The following last result concludes that, under Assumption~H, the
sequence of random variables $({\cal X}^N)_N$ valued in $\mathcal{M}_1(E)$ 
converges to the minimal quasi-stationary distribution of $X$. The
constant $\gamma$ and the measure $\nu_{QSD}$ are obtained from
Theorem~\ref{thm:main-QSD}.

\begin{thm}
\label{thm:main-cv}
Assume that Assumption~H is satisfied. Then there exists a constant
$d>0$ such that, for all bounded measurable function
$f:E\rightarrow\R$ and all
$N> \frac{\lambda_1}{\lambda_1-\|\kappa\|_\infty}$,
\[
  \E\left|{\cal X }^N(f)-\nu_{QSD}(f)\right|\leq \frac{d}{N^{\alpha}}\,\|f\|_\infty,
\]
where $\alpha=\frac{\gamma}{2(\|\kappa\|_\infty+\gamma)}$.

Moreover, for all function $f:E\rightarrow\R$ such that
$f(x)/V(x)\rightarrow 0$ when $V(x)\to +\infty$, we have
\[
  \E\left|{\cal X }^N(f)-\nu_{QSD}(f)\right|\xrightarrow[N\rightarrow+\infty]{} 0.
\]
\end{thm}

\section{Examples}
\label{sec:examples}

In this section, we apply our results to multi-dimensional birth and
death processes, continuous time Galton-Watson processes, continuous
time multi-di\-men\-sio\-nal Galton-Watson processes and multi-dimensional
diffusion processes.

\subsection{Multi-dimensional birth and death processes}
Let $Y$ be a continuous-time multitype birth and death process,
taking values in $E\cup\{\d\}=\mathbb N^d$ for some $d\geq 1$, with
transition rates
\[
  q_{x,y}=
  \begin{cases}
    b_i(x) & \text{if }y=x+e_i,\\ 
    d_i(x) & \text{if }y=x-e_i,\\
    0 & \text{otherwise,}
  \end{cases}
\]
where $(e_1,\ldots,e_d)$ is the canonical basis of $\mathbb Z_{+}^d$
(where $\mathbb Z_+=\{0,1,\ldots\})$), and $\d=(0,\ldots,0)$.  We
assume that $b_i(x)>0$ for all $1\leq i\leq d$ and
all $x\in E$, \textcolor{black}{and $d_i(x)>0$ for all $x\in E$ and $1\leq i\leq d$ such that $x_i\geq 1$ (of course, $d_i(x)=0$ otherwise).}

The following result provides a general explicit criterion on the
parameters of the process ensuring that the results of
Section~\ref{sec:mainresults} hold true.

\begin{prop}
  If
  \begin{align}
    \label{eq:multi-bd-1}
    \frac{1}{|x|}\sum_{i=1}^d\big(d_i(x)-b_i(x)\big) \to +\infty \quad\text{ when }|x|\to+\infty,
  \end{align}
  or if there exists $\delta>1$ such that
  \begin{align}
    \label{eq:multi-bd-2}
    \sum_{i=1}^d\big(d_i(x)-\delta\,b_i(x)\big)\to+\infty\quad\text{ when }|x|\to+\infty,
  \end{align}
  then Assumption~H is satisfied.
\end{prop}

Indeed, one can choose
$V(x)=|x|=x_1+\ldots+x_d$ if~\eqref{eq:multi-bd-1} is satisfied, and
$V(x)=\exp(\varepsilon x_1+\cdots+\varepsilon x_d)$ with
$\varepsilon>0$ small enough if~\eqref{eq:multi-bd-2} is satisfied
(see~\cite[Example~7]{ChampagnatVillemonais2017b} for the details)
and the fact that the killing rate is bounded by
$d_1(e_1)+\cdots+d_d(e_d)$. The rest of Assumption~H is a simple
consequence of the irreducibility of the process and the fact that
the state space is discrete.

Theorem~\ref{thm:main-QSD} was already obtained in~\cite{ChampagnatVillemonais2017b}, while, as
far we know, the convergence of $({\cal X}^N)_{N}$ for
multi-dimensional birth and death processes toward $\nu_{QSD}$ is
completely new.

\begin{rem}
The one-dimensional birth and death models considered in
Examples~2.7, 2.8 and~2.9 of~\cite[Section~2]{Villemonais2015} all
satisfy Assumption~H. These birth and death processes studied in the
above mentioned paper evolve in $\N=\{1,2,\ldots\}$, with positive
birth rates denoted by $b_x$, $x\in\N$, death rates denoted by $d_x$,
$x\geq 2$, and with $\kappa(x)=d_1\11_{x=1}>0$. The examples
considered are
\begin{itemize}
\item $b_x=b\,x^a$ and $d_x=d\,x^a$ for all $x\geq 1$, where $b<d$ are two positive constants and $a>0$ is fixed,
\item $b_i=b>0$  and $d_i=d>0$ for all $i\geq 2$, where $b<d$ are positive constants, with $d_1>0$ such that $(\sqrt{d}-\sqrt{b})^2> d_1$ and $b_1>0$,
\item $b_i=|\sin(i\pi/2)|i+1$ and $d_i=4i$ for all $i\geq 1$.
\end{itemize} 
The existence of the Lyapunov function $V$ is already proved in this
reference, while the rest of Assumption~H is a simple consequence of
the irreducibility of birth and death processes with positive
coefficients.

While we do not improve the extent of the domain of attraction of the
minimal quasi-stationary distribution, we emphasise that,
when~\cite{Villemonais2015} relies on the spectral theoretical results
derived in~\cite{DoornErik1991} and leads to weak convergence of
measures, our approach is probabilistic and we prove a stronger
exponential convergence toward the quasi-stationary distribution.

Theorem~\ref{thm:main-ergodicity} and the convergence in law of
$({\cal X}^N)_{N}$ toward $\nu_{QSD}$ (in the weak topology)  was
already proved in the above reference, but the speed of convergence of
Theorem~\ref{thm:main-cv} is new.
\end{rem}

\subsection{Continuous-time Galton Watson processes}
  A continuous time Galton-Watson process describes the evolution of a
  population where individuals reproduce independently at rate $1$ and whose
  progeny follows a common law $p$ on $\mathbb Z_+$ (here, the progeny
  replaces the parent, who dies during the reproduction event). More
  formally, its dynamic is described by the following  infinitesimal
  generator (with $\d=0$ the unique absorbing point and $f$ is any bounded function):
  \[
    L f(x)=\sum_{n\geq -1} q(x,x+n)(f(x+n)-f(x)),\ \forall x\in \Z_+
  \]
  where the jump rates matrix $q(\cdot,\cdot)$ is given by $q(0,\cdot)=0$ and
  \[
    q(x,x+n)= x\,p(n+1)\ \text{ for all $x\geq 1$ and $n\geq -1$}.
  \]

  Our first result concerns the existence of a minimal
  quasi-stationary distribution for this process.

  \begin{prop}
    Assume that $p(0)>0$ and $p(\{2,3,\ldots\})>0$, that $p$ admits a
    first moment $m=\sum_{n=0}^\infty n\,p(n)\in(0,1)$ and that there
    exists $\alpha>1$ such that
    $\sum_{n=0}^\infty n^\alpha\,p(n)<+\infty$. Then
% the assumptions of~\cite[Theorem~5.1]{ChampagnatVillemonais2017b} are satisfied and 
    the continuous time Galton-Watson
    process with reproduction law $p$ admits a unique quasi-stationary
    distribution $\nu_{QSD}$ such that
    $\sum_{x=1}^\infty \nu_{QSD}(\{x\})\,x^\alpha<+\infty$. Moreover,
    there exist two positive constants $c,\gamma>0$ such that
    \begin{align*}
      \left\|\P_{\mu}(Y_t\in\cdot\mid t<\tau_\d)-\nu_{QSD}\right\|_{TV}\leq c\,\mu(V)\,e^{-\gamma t},\ \forall t>0,
    \end{align*}
    for any probability distribution $\mu$ on $E$ such that
    $\sum_{x=1}^\infty \mu(\{x\})\,x^\alpha<+\infty$.
  \end{prop}

  A similar result was obtained in~\cite{AsselahFerrariEtAl2016} using
  different methods with the requirement that $\alpha=2$. Our result
  provides a sharper result by allowing any reproduction law with a
  moment of order $\alpha>1$.
% , where $\alpha$ is chosen arbitrary small in $(1,+\infty)$.
  Our result also provides the additional exponential
  convergence in total variation norm and all the consequences listed
  in~\cite{ChampagnatVillemonais2017b}.

  Before proving the above proposition, we state a criterion implying
  that the continuous time Galton-Watson process satisfies
  Assumption~H.

  \begin{prop}
    Assume that $p(0)>0$ and $p(\{2,3,\ldots\})>0$, that $p$ admits a
    first moment $m=\sum_{n=0}^\infty n\,p(n)\in(0,1)$ and that there
    exists $\alpha>p(0)/(1-m)$ such that
    $\sum_{n=0}^\infty n^\alpha\,p(n)<+\infty$. Then Assumption~H
    holds true.
  \end{prop}

\textcolor{black}{Note that, since $m>1-p(0)$, the assumption on $\alpha$ implies that $\alpha>1$.}
  Thanks to Theorem 1.1 in~\cite{AsselahFerrariEtAl2016}, the ergodicity of
  the Fleming-Viot process and its weak convergence toward the
  quasi-stationary distribution were already known with the condition
  that $p$ admits an exponential moment. It is thus improved here by
  considering reproduction laws admitting a polynomial moment of
  explicit order and by providing convergence in a stronger sense. One
  major advantage of our approach is also its flexibility : we will
  show in the next section how this result easily generalises to
  continuous time multi-type Galton-Watson processes.

  \medskip We first prove the second proposition, because its proof is
  straightforward. Setting $V(x)=x^\alpha$, one obtains for all
  $x\geq 1$,
  \[
    LV(x)= V(x) \sum_{n=0}^\infty x\,p(n)\,\left[\left(1+\frac{n-1}{x}\right)^\alpha-1\right]\sim_{x\to+\infty} \alpha (m-1)\,V(x).
  \]
  Since we assumed that $\alpha(m-1)< -p(0)$ and since $p(0)$ is the maximum of the
  absorption rate in this model, one deduces that the Lyapunov
  Assumption~H is satisfied. As in the previous example, the rest of Assumption~H is a simple
  consequence of the irreducibility of the process and the fact that
  the state space is discrete.

  Let us now prove the first proposition. \textcolor{black}{This is a direct consequence of Theorem~5.1
    in~\cite{ChampagnatVillemonais2017b}. To apply this result, we only need} to prove that the constant
  \[
    \lambda_0:=\inf\left\{\lambda>0,\text{ s.t. }\liminf_{t\rightarrow+\infty}e^{\lambda t}\,\P_{x}\left(Y_t=x\right)>0\right\}
  \]
  is strictly smaller than $\alpha(1-m)$.
%  (where $\alpha>1$ under our
%  assumptions) (indeed, once this is proved, the proposition is an
%  immediate consequence of).
  Note that,
  because of the irreducibility of the process and the discreteness of
  the state space, this constant does not depend on $x\in \N$ and we also have for all $L\in\N$
%   , because of the finiteness of the sets $\{1,\ldots,L\}$ where $L\in\N$, we also have
  \[
    \lambda_0=\inf\left\{\lambda>0,\text{ s.t. }\exists L>0,\,\liminf_{t\rightarrow+\infty}e^{\lambda t}\,\P_{x}\left(Y_t\in\{1,\ldots, L\}\right)>0\right\}
  \]
  % In order to prove that $\lambda_0< \alpha(1-m)$, fix $\lambda\in(1-m,\alpha(1-m))$
  % and let us prove that there exists $L>0$ such that
  % \[
  %   \liminf_{t\rightarrow+\infty}e^{\lambda t}\,\P_{x}\left(Y_t=x\right)>0.
  % \]

  We first observe that, setting $\varphi(x)=x$ for all $x\in \N$, we
  have $L \varphi(x)=(m-1)\varphi(x)$. Dynkin's formula implies that
  the process $M$ defined by $M_t=e^{(1-m)t} \varphi(Y_t)$ is a local
  martingale. Now, using the fact that, for all $t\geq 0$,
  $\E(V(Y_s))$ is uniformly bounded over $s\in[0,t]$ (because $LV$ is
  uniformly bounded from above) and the fact that $\varphi=o(V)$, one
  deduces that $M$ is uniformly integrable over $[0,t]$. In particular
  it is a martingale and one deduces that
  \begin{equation}
    \label{eq:mart}
    \E_x(\varphi(Y_t))=e^{(m-1)t}\varphi(x)\text{ for all $x\in \N$ and all $t\geq 0$.}
  \end{equation}

  % Now, straightforward calculations show that, for all $\varepsilon>0$ such
  % that $1+2\varepsilon< \alpha$ and setting
  % $V_\varepsilon(x)=x^{1+\varepsilon}$, one has
  % \begin{align*}
  %   LV_\varepsilon(x)&\leq \sum_{n=0}^{+\infty}p(n)\left[n-1+\varepsilon x \ln\left(1+\frac{n-1}{x}\right)\left(1+\frac{n-1}{x}\right)^{1+\varepsilon}\right]\,V_\varepsilon(x)\\
  %                    &\leq \left[m-1+\varepsilon\, C_{st}\,x\,\sum_{n=0}^{+\infty}p(n)\left(1+\frac{n-1}{x}\right)^{\alpha}\right]\,V_\varepsilon(x),
  % \end{align*}
  % where $C_{st}$ is a constant (which, in particular, does not depend
  % on $\varepsilon$). Now, one easily checks that
  % $\sum_{n=0}^{+\infty}p(n)\left(1+\frac{n-1}{x}\right)^{\alpha}$ is
  % uniformly bounded in $x\in\N$ and hence there exists $\varepsilon>0$ such that
  % \[
  %   LV_{\varepsilon}(x)\leq -\lambda V_{\varepsilon}(x),\text{ for all $x\in \N$.} 
  % \]
  % Similarly as for the computation of $\E(\varphi(Y_t))$, one deduces
  % that
  % \[
  %   \E_x(V_\varepsilon(Y_t))\leq e^{-\lambda t}
  %   V_\varepsilon(x),\text{ for all $x\in \N$.}
  % \]
  % We hence have, for all $L\geq 1$,
  % \begin{align*}
  %   e^{\lambda t} \E_x(\varphi(Y_t)) & = e^{\lambda t} \E_x(\varphi(Y_t)\11_{V(Y_t)\leq L \varphi(Y_t)})+e^{\lambda t}  \E_x(\varphi(Y_t)\11_{V(Y_t)\geq L \varphi(Y_t)})\\
  %                                    & \leq e^{\lambda t} \E_x(\varphi(Y_t)\11_{V(Y_t)\leq L \varphi(Y_t)}) + e^{\lambda t} \E_x(\frac{V(Y_t)}{L})
  % \end{align*}
  
  Now, using the fact that
  $LV(x)\leq \alpha(m-1) V(x)+C \varphi(x)$ for some positive
  constant $C>0$,
% that may change from line to line, 
  one deduces, using Dynkin's formula  that
\textcolor{black}{
  \begin{align*}
    \E_x\left(e^{(1-m) t}V(Y_t)\right) & \leq V(x) + \int_0^t \left[(\alpha-1)(m-1) \E_x(e^{(1-m)s}V(Y_s)) +
      C\,\E_x(e^{(1-m)s}\varphi(Y_s))\right] \,ds \\
    & \leq (1+C)V(x) - (\alpha-1)(1-m)\int_0^t \E_x(e^{(1-m)s}V(Y_s)) \,ds.
  \end{align*}
  Reminding that $\alpha>1$, we deduce that}
  % and hence, using classical Grownall type inequalities (see for instance
  % Lemma~7 in~\cite{MaillerVillemonais2018}) that
  \[
     \E_x\left(e^{(1-m) t}V(Y_t)\right) \leq (1+C) V(x),\ \forall t\geq 0.
  \]
  From this inequality and~\eqref{eq:mart}, one obtains, for all \textcolor{black}{$L\in\N$ and all $t\geq 0$,
  \begin{align*}
    \varphi(x)&=e^{(1-m) t}\E_x(\varphi(Y_t))\\
              &\leq e^{(1-m) t}\max_{y\in \{1,\ldots,L\}}\varphi(y)\,\P_x(Y_t\in \{1,\ldots,L\}) +e^{(1-m) t} \E_x\left(\frac{V(Y_t)}{L^{\alpha-1}}\right)\\
              &\leq e^{(1-m) t}\max_{y\in \{1,\ldots,L\}}\varphi(y)\,\P_x(Y_t\in \{1,\ldots,L\}) + \frac{1+C}{L^{\alpha-1}}\, V(x).
  \end{align*}
  For a fixed $x\in \N$, one can choose $L$ large enough so that
  $(1+C)V(x)/L^{\alpha-1}\leq \varphi(x)/2$ and deduce} that
  \[
    \liminf_{t\rightarrow+\infty} e^{(1-m) t} \P_x(Y_t\in \{1,\ldots,L\}) >0,
  \]
  which implies that $\lambda_0 \leq (1-m)$ and concludes the proof.

  \begin{rem}
    Note that equation~\eqref{eq:mart} also immediately implies that
    $\lambda_0\geq (1-m)$ so that we in fact proved that
    $\lambda_0=1-m$. The approach employed here for the study of
    quasi-stationary distributions is of course a general strategy
    relying on~\cite{ChampagnatVillemonais2017b} and on Dynkin's formula that may be useful for
    the study of several processes. The main point is to dispose of a
    super harmonic function $\varphi$ and a Lyapunov-type norm-like
    function $V$ that dominates $\varphi$ at infinity.
  \end{rem}

  \subsection{Multi-type continuous time Galton-Watson processes}

  % A nice feature of our result is that it transfers very well by
  % domination arguments. The aim of this section is to show how this
  % principle applies for multi-type Galton-Watson processes.

  A multi-type continuous time Galton-Watson process describes the
  evolution of a population of typed individuals, with $d\geq 1$
  possible types denoted by $1,\ldots,d$. These individuals reproduce
  independently at rate $\lambda_i>0$ and their progeny has a law $p_i$ on
  ${\mathbb Z}^d_+$ (here again, the progeny replaces the parent, who
  dies during the reproduction event), depending on the type
  $i\in\{1,\ldots,d\}$ of the parent. Its dynamic is thus described by
  the following infinitesimal generator (where
  $f$ is any bounded function):
  \[
    L f(x)=\sum_{i=1}^d\sum_{n\in\Z_+^d} q_i(x,x+ n-e_i)(f(x+ n-e_i)-f(x)),\ \forall x\in \Z^d_+
  \]
  where $(e_1,\ldots,e_d)$ is the canonical basis of $\Z^d_+$ and the
  jump rates matrix $q_i(\cdot,\cdot)$ is given by $q_i(x,\cdot)=0$ if
  $x_i=0$ and
  \[
    q_i(x,x+n-e_i)= \lambda_i\, x_i\,p_i(n)\text{ otherwise.}
  \]
  \textcolor{black}{We also assume that $p_i$ is such that the associated Galton-Watson process is irreducible in
    $\Z_+^d\setminus\{0\}$, so that $\d=(0,\ldots,0)$ is the unique absorbing subset of $\Z_+^d$.}

  % In the following proposition, we say that a $d$-type reproduction
  % law $(p_i)_{i\in\{1,\ldots,d\}}$ is dominated by a $1$-type law $p$ on
  % $\Z_+$ if and only if
  % \[
  %   \sum_{n\in\Z_+^d,\ |n|\geq x} p_i(n)\leq \sum_{y\geq x} p(y),\ \forall i\in\{1,\ldots,n\},\ \forall x\in\Z_+,
  % \]
  % where $|\cdot|$ is the $l_1$ norm.  This means that the size of the
  % progeny of one individual of any type $i\in\{1,\ldots,d\}$ is
  % stochastically dominated by the law $p$.

  In the following proposition, given such a rate $\lambda_i$ and a
  law $p_i$ on $\Z_+^d$, we denote by $M$ the matrix of mean
  offsprings defined y
  $M_{ij} =\lambda_i\,\sum_{n\in\Z_+^d} n_j p_i(n)$ and by $Q$ the
  matrix defined by $Q_{ij}=M_{ij}-\delta_{ij} \lambda_i$. Note that
  the irreducibility of the Galton-Watson process implies the
  irreducibility and aperiodicity of $e^{tQ}$ and hence, by
  Perron-Frobenius theorem, the existence of a positive
  right-eigenvector $v\in(0,+\infty)^d$ associated to the spectral
  radius $\rho\in\R$. Recall that $\rho<0$ corresponds to the
  subcritical case.

  % the matrix $Q$ defined by $Q_{ij}=M_{ij}-\delta_{ij} \lambda_i$ is
  % irreducible (since it is assumed that the Galton-Watson process is
  % irreducible). We recall that
  % $M_{ij} =\lambda_i\,\sum_{n\in\Z_+^d} n_j p_i(n)$ and that, by
  % Perron-Frobenius theorem, the irreducibility of $Q$ implies the
  % existence of a positive right-eigenvector $v\in(0,+\infty)^d$
  % associated to the spectral radius $\rho>0$.

  % the
  % moment of order $\alpha$ of $p_i$ is defined as
  % $\sum_{n\in\Z_+^d} |n|^\alpha\,p_i(n)$, where $|\cdot|$ is the $l_1$
  % norm.
  
  \begin{prop}
    Assume that $(\lambda_i)_{i\in\{1,\ldots,d\}}$ and
    $(p_i)_{i\in\{1,\ldots,d\}}$ are such that the associated
    multi-type continuous time Galton-Watson process is irreducible in
    $\Z_+^d\setminus \{0\}$. Assume also that the spectral radius
    $\rho$ of the matrix $Q$ is negative and that $p_i$ admits a
    moment of order $\alpha> \max_i p_i(0)/(-\rho)$, then Assumption~H
    is satisfied by the multi-type continuous time Galton-Watson
    process with reproduction laws $(p_i)_{i\in\{1,\ldots,d\}}$ and
    reproduction rates $(\lambda_i)_{i\in\{1,\ldots,d\}}$.
  \end{prop}
  
  Since the state space is discrete and the process is assumed to be
  irreducible, the only difficulty is to prove the existence of a
  function $V$ such that $LV(x)\leq -\lambda_1 V(x)+C$ for some
  constants $\lambda_1>\max_i p_i(0)$ (note that $\max_i p_i(0)$ is
  the maximal rate of absorption for our process). In order to do so,
  one simply observes that, setting $V(x)=\left(\sum_{i=1}^d v_i x_i\right)^\alpha$ for all
  $x\in\Z_+^d$, 
  \begin{align*}
    LV(x)&= V(x)\,\sum_{i=1}^d \sum_{n\in\Z_+^d} x_i\,\lambda_i\,p_i(n)\,\left[\left(1+\frac{\sum_{j=1}^d v_j (n_j-\delta_{ij})}{\sum_{j=1}^d v_j x_j}\right)^\alpha-1\right]\\
         &\sim_{|x|\to+\infty} \frac{\alpha\,V(x)}{\sum_{j=1}^d v_j x_j}\,\sum_{i=1}^d x_i\,\sum_{j=1}^d\,\sum_{n\in\Z_+^d} p_i(n)\, v_j \lambda_i\,(n_j-\delta_{ij}) \\
         &= \frac{\alpha\,V(x)}{\sum_{j=1}^d v_j x_j}\,\sum_{i=1}^d x_i \sum_{j=1}^d q_{ij} v_j =  \alpha \rho V(x),
  \end{align*}
  so that $V$ satisfies the Lyapunov-type inequality of Assumption~H. 
  
% \textcolor{red}{Est-ce que {\c c}a ne vaudrait pas le coup de pr\'esenter un r\'esultat avec moments exponentiels parce que dans ce
%   cas, il y a des chances qu'on ait besoin d'hypoth\`eses moins restrictives sur la matrice des moyennes de la loi de reproduction ?}
% \textcolor{blue}{A discuter!}

  \subsection{Diffusion processes with soft killing}

  In this section, we consider the case of a multi-dimensional
  diffusion process with soft killing and provide a sufficient
  criterion for Assumption~H to apply.  Let $(X_t)_{t\in[0,+\infty)}$
  be the solution in $E=\R^d$ to the stochastic differential equation
  \begin{align*}
     dX_t= \sigma(X_t) \,dB_t+b(X_t)  dt,
  \end{align*}
  where $B$ is a standard $r$-dimensional Brownian motion,
  $\sigma:\R^d\rightarrow \R^{d\times r}$ and $b:\R^d\mapsto \R^d$ are
  locally H\"older-continuous and $a=\sigma\sigma^*$ is bounded and locally
  uniformly elliptic. The process $Y$ is subject to an additional soft
  killing $\kappa:\R^d\mapsto [0,+\infty)$, which is assumed to be
  uniformly bounded.

  \begin{prop}
    If there exist $\beta>0$ and $\gamma>0$ such that
    \begin{align*}
      \limsup_{|x|\rightarrow+\infty} \frac{\langle b(x),x\rangle}{\langle x,x\rangle^{1/2}}
      \leq -\beta\quad\text{ and }\quad \sum_{ij} a_{ij}(x) x_i x_j\leq \gamma \langle x,x\rangle,
      %-\frac{3}{2} \|\kappa\|^{1/2}_\infty,
    \end{align*}
    and such that $\beta^2>2\gamma\|\kappa\|_\infty$, then
    Assumption~H is satisfied with
    $V:x\in\R^d\mapsto \exp(\rho \langle x,x\rangle^{1/2})$, where
    $\rho>0$ is such that
    $\rho^2\gamma/2+\|\kappa\|_\infty < \beta\rho$.
  \end{prop}

  The fact that Assumption~H is satisfied is in fact a consequence of
  the material included in Section~4.2, Example~3
  of~\cite{ChampagnatVillemonais2017b}: the property
  on small sets and absorption probabilities are both consequences of
  the Harnack inequalities proved therein (where they are used to
  prove the similar properties (F1) and (F3)). The computation of $LV$ is only
  provided in a simpler case, so we give details below: fix $\varepsilon>0$ such that
    $\rho^2\gamma/2+\|\kappa\|_\infty < (\beta-2\varepsilon)\rho$. Then, for $|x|$ large enough,
  \begin{align*}
    L V(x) & =\sum_{i=1}^d \rho V(x)\frac{\langle b(x),x\rangle}{\langle x,x\rangle^{1/2}}+\frac{1}{2}\sum_{1\leq i,j\leq d}
    a_{ij}(x) V(x)\left[\rho^2\frac{x_i x_j}{\langle x,x\rangle}-\rho\frac{x_i x_j}{\langle
        x,x\rangle^{3/2}}+\rho\frac{\delta_{ij}}{\langle x,x\rangle^{1/2}}\right] \\ 
    & \leq V(x)\left( (-\beta+\varepsilon)\rho+\frac{\rho^2\gamma}{2}+\frac{d\rho{}\|a\|_\infty}{\langle x,x\rangle^{1/2}}\right).
  \end{align*}
  Therefore, $LV(x)\leq -(\|\kappa\|_\infty+\varepsilon\rho) V(x)$ for $|x|$ large enough, so $V$
  satisfies the Lyapunov-type part of Assumption~H.

\section{Proofs}
\label{sec:proofs}

\subsection{Proof of Theorem~\ref{thm:main-QSD}}

\textit{Step 1. Quasi-stationary behaviour of $(Y_n)_{n\in\N}$.}  Let
us first prove that the skeleton $(Y_n)_{n\in\N}$ admits a minimal
quasi-stationary distribution, by proving that Assumption~E
in~\cite{ChampagnatVillemonais2017b} is satisfied and hence that
Theorem~2.1 therein applies. Fix
$\lambda\in(\lambda_1,\|\kappa\|_\infty)$. One deduces from
Assumption~H, Dynkin's formula, classical upper bound techniques (see
for instance Lemma~7 in~\cite{MaillerVillemonais2018}) that, for all $x\in E$ and all
$t\geq 0$,
\begin{align}
  \label{eq:dynkin}
  \E_x(V(Y_t)\11_{t<\tau_\d})\leq \E_x(V(X_t)) \leq
  \frac{C}{\lambda_1-\lambda} \vee (e^{-\lambda t} V(x)).
\end{align}
We define the set $K:=\{x\in E,\ V(x)\leq (2\vee e^{\lambda})C/(\lambda_1-\lambda)\}$,
which is relatively compact by Assumption~H. Our aim is to prove Assumption~E
in~\cite{ChampagnatVillemonais2017b} with this set $K$, $\varphi_1=V$, $\varphi_2=\11_E$,
$\theta_1=e^{-\lambda}$ and $\theta_2=e^{-\|\kappa\|_\infty}$.

Let us first check that Assumption~(E1) holds for any
$n_1\geq t_K+\frac{\ln 2}{\lambda}$ (here $t_K>0$ is the
positive time appearing in the small set property for $K$). \textcolor{black}{Note that
  Assumption~(E1) is actually simply the small set property for $K$, but we will need this stronger version in the sequel.} 
We deduce
from the definition of $K$ and from the above inequality that, for all
$x\in K$,
\[
  \P_x(X_{n_1-t_K}\notin K)\leq  \frac{\E_x(V(X_{n_1-t_K}))}{2C/(\lambda_1-\lambda)}\leq 1/2,
\]
In particular, using Markov's property at time $n_1-t_K$ and the fact
that $K$ is a small set for $X$, we obtain
\[
  \P_x(X_{n_1}\in A)\geq \frac{\alpha_K}{2}\nu_K(A\cap K),
\]
and hence that
\begin{align}
  \label{eq:E1}
  \P_x(Y_{n_1}\in A)\geq \frac{\alpha_Ke^{-\|\kappa\|_\infty
  n_1}}{2}\nu_K(A\cap K)
\end{align}
for all $n_1\geq t_K+\frac{\ln 2}{\lambda}$.
% which entails (E1).

We now prove that (E2) is satisfied. Using~\eqref{eq:dynkin}, we
obtain, for all $x\in E$,
\[
  \E_x(V(Y_1)\11_{1<\tau_\d})\leq \theta_1 V(x)+\frac{2C}{\lambda_1-\lambda}\11_{K}(x).
\]
Also, we have
\[
  \E_x(\11_{Y_1\in E}\11_{1<\tau_\d})\geq e^{-\|\kappa\|_\infty}\E_x(\11_{X_1\in E})=e^{-\|\kappa\|_\infty}.
\]
Hence, since we set $\varphi_1=V$, $\varphi_2=\11_E$,
$\theta_1=e^{-\lambda}$ and $\theta_2=e^{-\|\kappa\|_\infty}$, one
deduces that Assumption~(E2) is satisfied.

Since $K$ is a relatively compact set, Assumption~(E3) is an immediate
consequence of the last part of Assumption~H.

Finally, setting $n_4(x)=\lceil t_K+\frac{\ln 2}{\lambda}\rceil$ for all
$x\in K$, one deduces from~\eqref{eq:E1} that Assumption~(E4) holds true.

We deduce from~\cite[Theorem~2.1]{ChampagnatVillemonais2017b} that there exist a
quasi-stationary distribution $\nu_{QSD}$ for $(Y_n)_{n\in\N}$ and two
constants $c'>0$ and $\gamma>0$ such that, for all probability measure
$\mu$ on $E$ and all $n\geq 0$,
\begin{align}
  \label{eq:QSD-step1}
  \left\|\P_\mu(Y_n\in\cdot\mid n<\tau_\d)-\nu_{QSD}\right\|_{TV}\leq c'\mu(V)e^{-\gamma n}.
\end{align}

\medskip\textit{Step 2. Conclusion.} Our aim is now to prove that
$(Y_t)_{t\geq 0}$ exhibits the same quasi-stationary behaviour as
$(Y_n)_{n\in\N}$. Let $\mu$ be probability measure on $E$ such that
$\mu(V)<\infty$, fix $t\geq 0$ and set $s=t-\lfloor t\rfloor\in[0,1)$. The law
$\mu_s$ of $Y_s$ satisfies (we use the fact that the killing rate of $Y$ is
bounded by $\|\kappa\|_\infty$, the fact that $LV\leq C$ and Dynkin's
formula)
\[
 \mu_s(V)= \E_\mu(V(Y_s)\mid s<\tau_\d)\leq e^{\|\kappa\|_\infty s} \E_\mu(V(X_s))\leq e^{\|\kappa\|_\infty +C} \mu(V). 
\]
From~\eqref{eq:QSD-step1}, we deduce that
\[
  \left\|\P_{\mu_s}(Y_{\lfloor t\rfloor}\in\cdot\mid \lfloor t\rfloor<\tau_\d)-\nu_{QSD}\right\|_{TV}\leq c'e^{\|\kappa\|_\infty +C}\mu(V)e^{-\gamma \lfloor t\rfloor}.
\]
Setting $c=c'e^{\|\kappa\|_\infty +C+\gamma}$, one deduces from Markov
property at time $s$ and from the above inequality that, for all
$t\geq 0$,
\[
  \left\|\P_\mu(Y_t\in\cdot\mid t<\tau_\d)-\nu_{QSD}\right\|_{TV}\leq c\mu(V)e^{-\gamma t}
\]
This concludes the proof of Theorem~\ref{thm:main-QSD}.

\begin{rem}
  One could have been tempted to use Theorem~3.5
  in~\cite{ChampagnatVillemonais2017b} which originally deals with
  continuous time Markov processes and allows bounded functions for
  the equivalent of $V$. However, this latitude comes with a far
  greater complexity and, in particular, one needs to check that the
  strong Markov property is satisfied at the entry time of $K$ (or, in
  the present case, a superset of $K$ which is also a small set),
  which can be quite challenging to prove in general (especially for
  discontinuous processes).

  Since we only consider norm-like functions $V$ in this paper, we can
  drop this technical requirement, which simplifies the verification
  of our assumptions. Note that the proofs of
  Theorems~\ref{thm:main-ergodicity} and~\ref{thm:main-cv} also make
  use of the fact that $V$ is a norm-like function.
\end{rem}

\subsection{Proof of Theorem~\ref{thm:main-ergodicity}}

We define the function $f:E^N \rightarrow \R$ by $f(x)= \sum_{i=1}^N V(x_i)$.
We have, for all $x=(x_1,\ldots,x_N)\in E^N$,
\begin{align*}
L^Nf(x)&=\sum_{i=1}^N L V(x_i)+\sum_{i=1}^N \kappa(x_i)\,\frac{1}{N-1}\sum_{j=1,\,j\neq i}^N V(x_j)-V(x_i).
\end{align*}
Assumption~H implies that
\begin{align*}
L^Nf(x)&\leq  -\left(\lambda_1-\frac{N}{N-1}\|\kappa\|_\infty\right)f(x)+C,
\end{align*}
where $\lambda_1-\frac{N}{N-1}\|\kappa\|_\infty$ is positive since
$N >\frac{\lambda_1}{\lambda_1-\|\kappa\|_\infty}$ by assumption.

\noindent For any constant $k>0$, the set of $N$-tuples $x\in E^N$
such that $f(x)\leq k$ is a small set for the Fleming-Viot process,
because the level sets of $V$ are small sets for $X$ and because of
the boundedness of $\kappa$.  Thus, using the Foster Lyapunov
criterion of \cite[Theorem 6.1, p.536]{MeynTweedie1993} (see
also~\cite[Proposition~1.4]{Hairer2010} for a simplified account on
the subject), we deduce that the Fleming-Viot process is exponentially
ergodic and, denoting by $M^N$ its unique stationary distribution, we
also have
\begin{align}
\label{eq:tightness}
\int_{E^N} f(x)\,dM^N(x)\leq C/\left(\lambda_1-\frac{N}{N-1}\|\kappa\|_\infty\right).
\end{align}
This concludes the proof of Theorem~\ref{thm:main-ergodicity}.

\subsection{Proof of Theorem~\ref{thm:main-cv}}

Since the rebirth rate $\kappa$ is uniformly bounded, it is well known (see for instance \cite{Villemonais2014,CloezThai2016,CerouDelyonEtAl2016})
that there exists a constant $d_0>0$ such that, for all $N\geq 2$ and all bounded measurable functions $f:E\rightarrow \R$, 
\[
  \E\left|\mu^N_t(f)-\E_{\mu_0^N}\left(f(Y_t)\mid t<\tau_\d\right)\right|\leq \frac{d_0}{\sqrt{N}} \|f\|_\infty e^{\|\kappa\|_\infty t}.
\]
In particular, one deduces that, for all $t\geq 0$ (we use the
stationarity of $M^N$ for the first and second equality and the results of
Theorem~\ref{thm:main-QSD} for the second inequality),
\begin{align*}
  \E\left(|{\cal X }^N(f)-\nu_{QSD}(f)|\right)&=\E_{M^N}\left|\mu^N_t(f)-\nu_{QSD}(f)\right|\\
                                              &\leq \E_{M^N}\left|\mu^N_t(f)-\E_{\mu_0^N}\left(f(Y_t)\mid t<\tau_\d\right)\right|+\E_{M^N}\left|\E_{\mu_0^N}\left(f(Y_t)\mid t<\tau_\d\right)-\nu_{QSD}(f)\right|\\
                                              &= \E_{M^N}\left|\mu^N_t(f)-\E_{\mu_0^N}\left(f(Y_t)\mid t<\tau_\d\right)\right|+\E\left|\E_{{\cal X}^N}\left(f(Y_t)\mid t<\tau_\d\right)-\nu_{QSD}(f)\right|\\ 
                                              &\leq \frac{d_0}{\sqrt{N}} \|f\|_\infty e^{\|\kappa\|_\infty t} + c \|f\|_\infty e^{-\gamma t} \E\left({\cal X}^N(V)\right)\\
                                              &\leq \frac{d_0}{\sqrt{N}} \|f\|_\infty e^{\|\kappa\|_\infty t} +C  c \|f\|_\infty e^{-\gamma t} /(\lambda_1-\|\kappa\|_\infty N/(N-1)),
\end{align*}
where the last inequality follows from
Theorem~\ref{thm:main-ergodicity}. In particular, choosing
$t=\frac{\ln N}{2(\|\kappa\|_\infty+\gamma)}$, one deduces that there
exists a constant $d>0$ such that
\[
  \E\left(|{\cal X }^N(f)-\nu_{QSD}(f)|\right)\leq d\|f\|_\infty/N^{\alpha},
\]
where $\alpha=\frac{\gamma}{2(\|\kappa\|_\infty+\gamma)}$.

In order to conclude the proof of Theorem~\ref{thm:main-cv}, fix
$f:E\rightarrow \R_+$ such that $f=o(V)$. Then, since
$\E({\cal X}^N(V))$ is uniformly bounded, one deduces that $f$ is
uniformly integrable with respect to this sequence of measures : for
all $\delta>0$, 
\[
  \E\left(|{\cal X }^N(f)-\nu_{QSD}(f)|\right)\leq \E\left(|{\cal X
    }^N(f\11_{V\leq \delta})-\nu_{QSD}(f\11_{V\leq \delta})|\right)+\E\left(|{\cal
      X }^N(f\11_{V> \delta})|\right)+\E\left(\nu_{QSD}(V\11_{V\leq \delta})\right)
\]
The first term in the right hand side goes to zero when
$N\rightarrow+\infty$ (because of the above analysis), while the
second and third terms converge to $0$ when $\delta\to +\infty$,
uniformly in $N$. This implies that the left hand side converges to
$0$ when $N\to +\infty$.

This concludes the proof of Theorem~\ref{thm:main-cv}.

\bibliographystyle{abbrv}
\bibliography{biblio-denis}

\begin{thebibliography}{10}

\bibitem{AsselahFerrariEtAl2011}
A.~Asselah, P.~A. Ferrari, and P.~Groisman.
\newblock Quasistationary distributions and {F}leming-{V}iot processes in
  finite spaces.
\newblock {\em J. Appl. Probab.}, 48(2):322--332, 2011.

\bibitem{AsselahFerrariEtAl2016}
A.~Asselah, P.~A. Ferrari, P.~Groisman, and M.~Jonckheere.
\newblock Fleming-{V}iot selects the minimal quasi-stationary distribution: the
  {G}alton-{W}atson case.
\newblock {\em Ann. Inst. Henri Poincar\'e Probab. Stat.}, 52(2):647--668,
  2016.

\bibitem{AsselahThai2012}
A.~{Asselah} and M.-N. {Thai}.
\newblock {A note on the rightmost particle in a Fleming-Viot process}.
\newblock {\em ArXiv e-prints}, Dec. 2012.

\bibitem{BieniekBurdzyEtAl2011}
M.~Bieniek, K.~Burdzy, and S.~Finch.
\newblock Non-extinction of a fleming-viot particle model.
\newblock {\em Probability Theory and Related Fields}, pages 1--40, 2011.

\bibitem{BurdzyHolystEtAl1996}
K.~Burdzy, R.~Holyst, D.~Ingerman, and P.~March.
\newblock Configurational transition in a fleming-viot-type model and
  probabilistic interpretation of laplacian eigenfunctions.
\newblock {\em J. Phys. A}, 29(29):2633--2642, 1996.

\bibitem{BurdzyHolystEtAl2000}
K.~Burdzy, R.~Ho{\l}yst, and P.~March.
\newblock A {F}leming-{V}iot particle representation of the {D}irichlet
  {L}aplacian.
\newblock {\em Comm. Math. Phys.}, 214(3):679--703, 2000.

\bibitem{CerouDelyonEtAl2016}
F.~{Cerou}, B.~{Delyon}, A.~{Guyader}, and M.~{Rousset}.
\newblock {A Central Limit Theorem for Fleming-Viot Particle Systems with Soft
  Killing}.
\newblock {\em ArXiv e-prints}, Nov. 2016.

\bibitem{ChampagnatVillemonais2017b}
N.~{Champagnat} and D.~{Villemonais}.
\newblock {General criteria for the study of quasi-stationarity}.
\newblock {\em ArXiv e-prints}, Dec. 2017.

\bibitem{ChampagnatVillemonais2017}
N.~{Champagnat} and D.~{Villemonais}.
\newblock {Lyapunov criteria for uniform convergence of conditional
  distributions of absorbed Markov processes}.
\newblock {\em ArXiv e-prints}, Apr. 2017.

\bibitem{CloezThai2016}
B.~Cloez and M.-N. Thai.
\newblock Quantitative results for the {F}leming-{V}iot particle system and
  quasi-stationary distributions in discrete space.
\newblock {\em Stochastic Process. Appl.}, 126(3):680--702, 2016.

\bibitem{DelyonCerouEtAl2017}
B.~{Delyon}, F.~{C{\'e}rou}, A.~{Guyader}, and M.~{Rousset}.
\newblock {A Central Limit Theorem for Fleming-Viot Particle Systems with Hard
  Killing}.
\newblock {\em ArXiv e-prints}, Sept. 2017.

\bibitem{FerrariMaric2007}
P.~A. Ferrari and N.~Mari{\'c}.
\newblock Quasi stationary distributions and {F}leming-{V}iot processes in
  countable spaces.
\newblock {\em Electron. J. Probab.}, 12:no. 24, 684--702 (electronic), 2007.

\bibitem{GrigorescuKang2004}
I.~Grigorescu and M.~Kang.
\newblock Hydrodynamic limit for a {F}leming-{V}iot type system.
\newblock {\em Stochastic Process. Appl.}, 110(1):111--143, 2004.

\bibitem{GrigorescuKang2012}
I.~Grigorescu and M.~Kang.
\newblock Immortal particle for a catalytic branching process.
\newblock {\em Probab. Theory Related Fields}, 153(1-2):333--361, 2012.

\bibitem{GroismanJonckheere2013a}
P.~{Groisman} and M.~{Jonckheere}.
\newblock {Front propagation and quasi-stationary distributions: the same
  selection principle?}
\newblock {\em ArXiv e-prints}, Apr. 2013.

\bibitem{GroismanJonckheere2013}
P.~Groisman and M.~Jonckheere.
\newblock Simulation of quasi-stationary distributions on countable spaces.
\newblock {\em Markov Process. Related Fields}, 19(3):521--542, 2013.

\bibitem{Hairer2010}
M.~Hairer.
\newblock Convergence of markov processes (lecture notes).
\newblock www.hairer.org/notes/Convergence.pdf, 2010.

\bibitem{MaillerVillemonais2018}
C.~{Mailler} and D.~{Villemonais}.
\newblock {Stochastic approximation on non-compact measure spaces and
  application to measure-valued P$\backslash$'olya processes}.
\newblock {\em ArXiv e-prints}, Sept. 2018.

\bibitem{MeynTweedie1993}
S.~P. Meyn and R.~L. Tweedie.
\newblock Stability of {M}arkovian processes. {III}. {F}oster-{L}yapunov
  criteria for continuous-time processes.
\newblock {\em Adv. in Appl. Probab.}, 25(3):518--548, 1993.

\bibitem{DoornErik1991}
E.~A. van Doorn.
\newblock Quasi-stationary distributions and convergence to quasi-stationarity
  of birth-death processes.
\newblock {\em Adv. Appl. Probab.}, 23(4):683--700, 1991.

\bibitem{Villemonais2011}
D.~Villemonais.
\newblock Interacting particle systems and yaglom limit approximation of
  diffusions with unbounded drift.
\newblock {\em Electronic Journal of Probability}, 16:1663--1692, 2011.

\bibitem{Villemonais2014}
D.~Villemonais.
\newblock General approximation method for the distribution of {M}arkov
  processes conditioned not to be killed.
\newblock {\em ESAIM Probab. Stat.}, 18:441--467, 2014.

\bibitem{Villemonais2015}
D.~Villemonais.
\newblock Minimal quasi-stationary distribution approximation for a birth and
  death process.
\newblock {\em Electron. J. Probab.}, 20:no. 30, 18, 2015.

\end{thebibliography}

\end{document}